\documentclass{amsart}
\usepackage{fdsymbol}

\usepackage{comment}

\DeclareFontFamily{OT1}{pzc}{}
\DeclareFontShape{OT1}{pzc}{m}{it}{<-> s * [1.100] pzcmi7t}{}
\DeclareMathAlphabet{\mathpzc}{OT1}{pzc}{m}{it}

\usepackage{tikz-cd}
\tikzcdset{arrow style=math font}

\usepackage{biblatex}
\usepackage{booktabs}
\bibliography{references}

    \usepackage{etoolbox}
    \patchcmd{\section}{\scshape}{\large\bfseries}{}{}
    \makeatletter
    \renewcommand{\@secnumfont}{\bfseries}
    \makeatother

\numberwithin{equation}{section}
\newtheorem{theorem}{Theorem}[section]
\newtheorem*{theorem*}{Theorem}

\theoremstyle{definition}
\newtheorem{question}{Question}

\newtheorem{remark}[theorem]{Remark}

\def\gg{\mathfrak{g}}
\def\QQ{\mathbb{Q}}
\def\ZZ{\mathbb{Z}}
\def\WW{\mathcal{W}}

\def\Bau{{\sf Bau}}

\def\epi{\twoheadrightarrow}

\def\CC{\mathpzc{C}}

\def\Gr{\mathsf{Gr}}

\title[An overview of rationalization theories]{An overview of rationalization theories\\ of non-simply connected spaces \\ and non-nilpotent groups}

\author{Sergei O. Ivanov} 
\address{
Laboratory of Modern Algebra and Applications,  St. Petersburg State University, 14th Line, 29b,
Saint Petersburg, 199178 Russia}
\email{ivanov.s.o.1986@gmail.com}

\thanks{This work is supported by the Ministry of Science and Higher Education of the Russian Federation, agreement 075-15-2019-1619.}

\begin{document}

\def\arraystretch{1.5}

\maketitle

\begin{abstract}
We give an overview of five rationalization theories for spaces (Bousfield--Kan's $\mathbb Q$-completion;  Sullivan's rationalization; Bousfield's homology rationalization;  Casacuberta--Peschke's $\Omega$-rationalization; G\'{o}mez-Tato--Halperin--Tanr\'{e}'s $\pi_1$-fiberwise rationalization) that extend the classical rationalization of simply connected spaces. We also give an overview of the corresponding rationalization theories for groups ($\mathbb Q$-completion; $H\mathbb Q$-localization; Baumslag rationalization) that extend the classical Malcev completion. 
\end{abstract}

\section*{Introduction} 
A simply connected space $X$ is called rational if one of two equivalent properties is satisfied: (1) homotopy groups $\pi_*(X)$ are  $\QQ$-vector spaces; (2) reduced integral homology groups $\tilde H_*(X,\ZZ)$ are $\QQ$-vector spaces. It is important for the rational homotopy theory that for a simply connected space $X$ there exists a universal map to a rational space $X\to X_\QQ$ that induces isomorphisms $\pi_*(X_\QQ)\cong \pi_*(X)\otimes \QQ$ and $\tilde H_*(X_\QQ,\ZZ)\cong \tilde H_*(X,\QQ).$  The space $X_\QQ$ is called rationalization of $X.$

This note is devoted to an overview of five rationalization theories of all spaces that extend the classical rationalization of simply connected spaces:
\begin{enumerate}
    \item  Bousfield-Kan's $\QQ$-completion $\QQ_\infty$ \cite{bousfield1972homotopy};
    \item Sullivan's rationalization ${\sf Sul}$ \cite{sullivan1977infinitesimal}, \cite{halperin2015rational}, \cite{bousfield1976pl}; 
    \item Bousfield's homology rationalization $L_{H\QQ}$ \cite{bousfield1975localization}; 
    \item Casacuberta---Peschke's $\Omega$-rationalization $L_{\Omega\QQ}$ \cite{casacuberta1993localizing};
    \item G\'{o}mez-Tato---Halperin---Tanr\'{e}'s $\pi_1$-fiberwise rationalization $L^{\pi_1}_{\QQ}$ \cite{gomez2000rational}.
\end{enumerate}
We do not aim to provide an overview of the history of the topic. We just want to point out that there are five different approaches to rationalization of a space list their properties and show some relations between them. 

For spaces of finite rational type the Sullivan rationalization coincides with the $\QQ$-completion
\begin{equation}
    {\sf Sul}(X)\cong \QQ_\infty X.
\end{equation}
However, all other pairs of functors are non-isomorphic even for spaces of finite rational type but there are natural transformations:
\begin{equation}L_\QQ^{\pi_1} \longrightarrow L_{\Omega\QQ} \longrightarrow L_{H\QQ} \longrightarrow \QQ_\infty.\end{equation}

Some properties of these constructions can be summarised in the following table. 

\begin{center}
\begin{tabular}{|c|c|c|c|c|c|}
\hline 
 &  $\QQ_\infty$ &  ${\sf Sul}$ & $L_{H\QQ}$  & $L_{\Omega\QQ}$ & $L^{\pi_1}_{\QQ}$  \\
\hline 
Idempotent &  $-$ & $-$ & $+$ & $+$ & $+$  \\
\hline
``Right'' for nilpotent spaces &  $+$ & $-$ &  $+$ & $+$ & $-$ \\ 
\hline
``Right'' for nilpotent spaces of finite type &  $+$ & $+$ &  $+$ & $+$ & $-$ \\ 
\hline  
 \begin{tabular}{@{}c@{}}Has a simplicial group description  \end{tabular}
&  $+$ & $-$ & $?$ & $+$ & $+$ \\ \hline 
Aspherical for the wedge of two circles & $+$ & $+$ & $?$ & $+$ & $+$ \\
\hline
\end{tabular}    
\end{center}

The three functors $L^{\pi_1}_{\QQ},$ $L_{\Omega\QQ},$ $L_{H\QQ}$ are localizations, and they are defined by  three classes of maps: \begin{equation}
\{\pi_1\text{-fiberwise rational eq.}\}\subseteq \{\Omega\text{-rational eq.} \}\subseteq \{\text{homology rational eq.}\},
\end{equation}
and three classes of spaces:
\begin{equation}
\{\pi_1\text{-fiberwise rational sp.}\}\supseteq \{\Omega\text{-rational sp.} \}\supseteq \{\text{homology rational sp.}\},    
\end{equation}
that we also describe. 

\ 

Since we deal with non-simply connected spaces, we focus attention on the fundamental groups of all these constructions. These rationalization theories for  spaces correspond to the following rationalization theories for groups: 
 $\QQ$-completion;
 $H\QQ$-localization;
Baumslag rationalization.

\begin{center}
\begin{tabular}{|c|c|}
\hline
Spaces & Groups  \\
\hline
\begin{tabular}{@{}c@{}} Bousfield-Kan's $\QQ$-completion \\ Sullivan's rationalization \end{tabular}  & $\QQ$-completion \\ 
\hline
Homology rationalization & $H\QQ$-localization \\ 
\hline
$\Omega$-rationalization & Baumslag rationalization\\
\hline 
$\pi_1$-fiberwise rationalization & --- \\
\hline 
\end{tabular}
\end{center}

\

These three rationalization theories for groups generalize the classical Malcev completion \cite{mal1949class}. The $\pi_1$-fiberwise rationalization of a space does not correspond to a rationalization theory on groups because it does not affect the fundamental group. All these constructions can be generalized to localizations in any set of primes $P$ but for simplicity we focus our attention on rationalizations.

\section*{Acknowledgements} I am grateful to Emmanuel Farjoun and Stephen Halperin for useful discussions.

\section{General localizations in categories} 
Let $\CC$ be a category. For a morphism $f:c_1\to c_2$ an object $l\in \CC$ is called $f$-local, if the map $f^*: \CC(c_2,l) \to \CC(c_1,l)$ is a bijection. An object is called $\WW$-local for some class of morphisms $\WW$, if it is $w$-local for any $w\in \WW.$
A $\WW$-localization of an object $c\in \CC$ is a morphism $w:c\to l$ from the class $\WW$ to a local object. It is easy to check that a $\WW$-localization satisfies two universal properties: (1) it is universal (initial) morphism to a $\WW$-local object with domain $c$; (2) it is universal (terminal) morphism from $\WW$ with domain $c.$ If the $\WW$-localization exists for any object in $\CC,$ it defines a coaugmented functor $L:\CC\to \CC$ i.e. a functor with a natural transformation $\eta:{\sf Id}\to L.$  Moreover, this coaugmented functor is a localization functor in the following sense: the maps $\eta L, L\eta:L\to L^2$ are equal isomorphisms. In particular, $L\cong L^2.$ In this case the full subcategory of $\WW$-local objects is a reflective subcategory and the restriction of the functor $L$ is the reflector. Moreover, this reflector induces an equivalence of categories
\begin{equation}
\CC[\WW^{-1}] \simeq {\sf Loc}(\WW),    
\end{equation}
where ${\sf Loc}(\WW)$ is the subcategory of local objects. 

Some of these definitions are also useful in the context of categories enriched over spaces (or simplicially enriched categories). We denote by ${\sf Map}_\CC(c_1,c_2)$ the hom-space between objects $c_1,c_2\in \CC.$  In the enriched setting we say that an object $l$ is $f$-local if $f^*:{\sf Map}_\CC(c_2,l) \to {\sf Map}_\CC(c_1,l)$ is a weak equivalence. The theory of $f$-localizations of spaces (in this enriched sense) can be found in \cite{farjoun2006cellular}, \cite{hirschhorn2009model}.

\section{The Malcev completion of nilpotent groups} 

\subsection{Definition} 
Recall that the \emph{lower central series} of a group $G$ is a series of normal subgroups given by 
\begin{equation}G_{n+1}=[G_n,G], \hspace{1cm} G=G_1 \supseteq G_2 \supseteq \dots. \end{equation}
We say that a group is $n$-nilpotent if $G_n=1.$ A group is called nilpotent of class $n$ if it is $n+1$-nilpotent but not $n$-nilpotent. A group is called \emph{nilpotent} if is is $n$-nilpotent for some $n.$

We say that a group $G$ is \emph{rational}  (or uniquely divisible, or complete) if for any $n\geq 1$ the map $G\to G, g\mapsto g^n$ is a bijection. For example, an abelian group is rational iff it a vector space over $\QQ.$ In particular, there is a natural way of rationalization of an abelian group $A$: the map \begin{equation}A\longrightarrow A\otimes \QQ, \hspace{1cm} a\mapsto a\otimes 1\end{equation}
is the universal map from an abelian group to a rational group. 

Malcev \cite{mal1949class}, \cite{mal1949nilpotent} (see also \cite{hilton1973localization}, \cite{warfield2006nilpotent}) extended the functor of rationalization from the category of abelian groups to the category of nilpotent groups
\begin{equation} -\otimes \QQ : (\text{nilpotent groups}) \longrightarrow (\text{nilpotent groups}).\end{equation}
This functor is called the Malcev completion. The Malcev completion  can be defined as the universal map to a rational nilpotent group
\begin{equation}G \longrightarrow G\otimes \QQ.\end{equation} 
The Malcev completion shares a lot of nice properties with the rationalization on abelian groups. For example, it takes short exact sequences to short exact sequences and central extensions to central extensions \cite[Prop.1.10]{hilton2016localization}. It induces an isomorphism on homology \cite[Th.4.8]{hilton1973localization}
\begin{equation}
H_*(G,\QQ) \cong H_*(G\otimes \QQ ,\ZZ),     
\end{equation}
and others. 
There is no such a nice generalization of this functor to all groups. 

\subsection{The Malcev correspondence}

The Malcev completion sends nilpotent groups to rational nilpotent groups. Rational nilpotent groups is a convenient class of groups because the category of rational nilpotent groups is equivalent to the category of nilpotent Lie algebras over $\QQ$ \cite[Th.9]{mal1949nilpotent}, \cite[Th.12.11]{warfield2006nilpotent}, \cite[Th.4.6]{baumslag2007lecture}, \cite[Appendix A.]{quillen1969rational}. 
\begin{equation}
(\text{rational nilpotent groups})\simeq (\text{nilpotent Lie algebras over }\QQ).
\end{equation}
Moreover, if $G$ is a rational nilpotent group and $\gg$ is the corresponding Lie algebra, then the associated graded Lie algebras ${\sf gr}(G)=\bigoplus_{n\geq 1} G_n/G_{n+1}$ and ${\sf gr}(\gg)=\bigoplus_{n\geq 1} \gg_n/\gg_{n+1}$ are isomorphic. It follows that the category of rational $n$-nilpotent groups is equivalent to the category of $n$-nilpotent Lie algebras
\begin{equation}\label{eq:ss}
(\text{rational $n$-nilpotent groups})\simeq (\text{$n$-nilpotent Lie algebras over }\QQ).    
\end{equation}
The functor from nilpotent Lie algebras to rational nilpotent groups is very explicit. For any nilpotent Lie algebra $\gg$ over $\QQ$ we define a new binary operation by the Baker–Campbell–Hausdorff formula
\begin{equation}\label{eq:BCH}
a*b= a+b+\frac{1}{2}[a,b]+\frac{1}{12}([a,[a,b]] - [b,[a,b]])+\dots     
\end{equation}
This sum is finite because the Lie algebra is nilpotent. The Lie algebra together with this new binary operation is the corresponding rational nilpotent group.

\section{$\QQ$-completion of groups} 
\subsection{Definitions}
$\QQ$-completion of a group $G$ can be defined as follows. For each $n$ we consider the nilpotent group $G/G_n$ and take its Malcev's completion $G/G_n\otimes \QQ.$ These groups form an inverse tower
\begin{equation}
 G/G_1 \otimes \QQ \longleftarrow G/G_2 \otimes \QQ \longleftarrow G/G_3 \otimes \QQ \longleftarrow   \dots.      
\end{equation}
The inverse limit of this tower is the $\QQ$-completion of $G$
\begin{equation} \widehat{G}_\QQ = \varprojlim (G/G_n \otimes \QQ).\end{equation} 

There are several equivalent definitions of the $\QQ$-completion. Instead of using the ordinary lower central series one can use the rational lower central series. It is defined as follows
\begin{equation}G_{n+1}^{\QQ}={\rm Ker}(G_n^{\QQ} \longrightarrow G_n^{\QQ}/[G_n^{\QQ},G] \otimes \QQ),\end{equation}
where $G_1^{\QQ}=G.$ The map $G/G_n\otimes \QQ\to G/G_n^\QQ \otimes \QQ$ is an isomorphism (it follows from the fact that $G_n^\QQ/G_n$ is torsion.) It follows that 
\begin{equation}\widehat{G}_\QQ = \varprojlim (G/G_n^\QQ \otimes \QQ).\end{equation}

Another equivalent definition is more categorical, it can be found in \cite[Ch.IV,\S 2.2]{bousfield1972homotopy}. Consider the category of all nilpotent rational groups ${\sf Nil}_\QQ$ and the category $G\!\downarrow\!{\sf Nil}_\QQ$ of homomorphisms $G\to N,$ where $N\in {\sf Nil}_\QQ.$ Then the $\QQ$-completion is the limit of the functor $G\!\downarrow\!{\sf Nil}_\QQ \longrightarrow {\sf Gr}$ given by $ (G\to N)\mapsto N$  
\begin{equation}
\widehat{G}_\QQ = \underset{G \to N\in  {\sf Nil}_\QQ }{\rm lim} N. 
\end{equation}

\subsection{(Non)-idempotency} Generally $\QQ$-completion is not idempotent. For example, if $F_\infty$ is the countably generated free group, then its double $\QQ$-completion is not isomorphic to its $\QQ$-completion. However the functor of completion is idempotent for the class of groups with finite dimensional $H_1(G,\QQ)$ \cite[\S 13]{bousfield1977homological}. Note that all finitely generated groups are in this class.  Moreover, the completion preserves this class of groups and for groups in this class there is an isomorphism $H_1(G,\QQ)\cong H_1(\widehat{G}_\QQ,\QQ).$

\subsection{$\QQ$-completion of the free group.} Let $F=F(X)$ be a free group generated by some set $X.$ Then its $\QQ$-completion can be described on the language of Lie algebras using  Baker–Campbell–Hausdorff formula. 

Denote by $L=L(X)$ the free Lie algebra over $\QQ$ generated by $X$. Then it has a natural gradding $L=\bigoplus_{i\geq 1} L^i,$ where $L^i$ is the vector space generated by commutators of weight $i.$ We consider the {\it complete free Lie algebra}, which is given by 
\begin{equation}
\widehat{L} = \prod_{i\geq 1} L^i.    
\end{equation}
Note that the lower central series of $L$ can be described as follows $L_n=\bigoplus_{i\geq n}L^i.$ So, the complete free Lie algebra can be described as 
\begin{equation}
\widehat{L}=\varprojlim L/L_n.    
\end{equation}
There is a natural structure of a group on this Lie algebra given by the  Baker-Campbell-Hausdorff formula \eqref{eq:BCH}. Then there is an isomorphism 
\begin{equation}
\widehat{F}_\QQ\cong (\widehat{L},*).    
\end{equation}
In order to prove this it is enough to prove that 
$F/F_n \otimes \QQ \cong (L/L_n, *).$
It follows from the Malcev correspondence because $F/F_n\otimes \QQ$ is the free object in the category of $n$-nilpotent rational groups and $L/L_n$ is the free object in the category of $n$-nilpotent Lie algebras over $\QQ.$

For a free group of rank $\geq 2$ the second homology group of the completion is non-trivial (and even uncountable) \cite[Th.1]{ivanov2019finite}
\begin{equation}
H_2(\widehat{F}_\QQ,\QQ )\ne 0.
\end{equation}
However nothing is known about the homology groups for $n\geq 3$

\begin{question} Is it true that
$H_n(\widehat{F}_\QQ,\QQ )= 0$ for a free group $F$ and $n\geq 3?$
\end{question}

\section{$H\QQ$-localization of groups}

\subsection{Definition} A homomorphism $f:G\to G'$ is called rationally $2$-acyclic (or $H\QQ$-homomorphism) if the map $H_1(G,\QQ)\to H_1(G',\QQ)$ is an isomorphism and the map $H_2(G,\QQ)\to H_2(G',\QQ)$ is an epimorphism. This definition is natural because a homomorphism $f:G\to G'$ is rationally $2$-acyclic if and only if there exists a map of spaces $F:X\to X'$ that induces an isomorphism of all homology groups $H_*(X,\QQ)\cong H_*(X',\QQ)$ such that $f:G\to G'$ is isomorphic to $\pi_1(F):\pi_1(X)\to \pi_1(X')$  \cite[Lemma 6.1]{bousfield1975localization}. Roughly speaking, the class of rationally $2$-acyclic morphisms is the image of the class of rational homology equivalences under $\pi_1.$

We denote by $\WW_{H\QQ}$ the class of rationally $2$-acyclic homomorphisms. Then a group is called $H\QQ$-local if it is $\WW_{H\QQ}$-local and $H\QQ$-localization is the $\WW_{H\QQ}$-localization. The $H\QQ$-localization exists for all groups \cite[Th.5.2]{bousfield1975localization} and defines a functor
\begin{equation}\label{eq:HQ-localization}
\ell_{H\QQ} : {\sf Gr} \longrightarrow {\sf Gr}.
\end{equation}

The class of $H\QQ$-local groups can be described more explicitly. We say that a central extension $E\epi G$ is rational if the kernel is a $\QQ$-vector space. Then the class of $H\QQ$-local groups is the least class containing the trivial group, closed under small limits and rational central extensions \cite[Th.3.10]{bousfield1977homological}. This description implies that an $H\QQ$-local group is rational. 

The functor $\ell_{H\QQ}$ is right exact \cite{akhtiamov2021right}. There is a transfinite limit construction of $H\QQ$-localization \cite{bousfield1977homological}.

\subsection{$H\ZZ$-localization of $G$-modules} There is a theory of $H\ZZ$-localizations of $G$-modules, which is similar to the theory of $H\QQ$-localizations of groups. Let $G$ be a fixed group and consider the category of modules ${\sf Mod}(G).$ Then a homomorphism $f:M\to M'$ is called {\it $1$-acyclic} if $H_0(G,M)\to H_0(G,M')$ is iso and $H_1(G,M)\to H_1(G,M')$ is epi. A module $M$ is called $H\ZZ$-local if it is local with respect to the class of $1$-acyclic homomorphisms. The class of $H\ZZ$-local modules can be described as the least class containing the trivial module, closed small limits and central extensions of modules (an extension of a $G$-module $E\epi M$ is called {\it central} if the action of $G$ on the kernel is trivial) \cite{bousfield1977homological}.

\section{Baumslag rationalization of groups} 

\subsection{Definition}
Recall that a group $G$ is called rational if the power map $G\to G,g\mapsto g^n$ is bijective for any $n\geq 1.$ It is easy to see that for a group $G$ the following statements are equivalent: (1) $G$ is rational; (2) $G$ is local with respect to the maps $\cdot n:\ZZ\to \ZZ$ for $n\geq 1;$ (3) $G$ is local with respect to the embedding $\ZZ\hookrightarrow \QQ. $

For any group $G$ there exists a universal map to a rational group that we call the Baumslag rationalization 
\begin{equation}
G \longrightarrow \Bau(G)  
\end{equation}
(see \cite{ribenboim1987equations}, \cite{peschke1987h}). This defines a coaugmented functor
\begin{equation}
\Bau : \Gr \longrightarrow \Gr.
\end{equation}

The functor of Baumslag rationalization satisfies the following properties:

\begin{itemize}
\item $\Bau$ sends epimorphisms to epimorphisms;
\item Moreover, $\Bau$ is right exact \cite{akhtiamov2021right};
\item $|\Bau(G)|\leq |G|$ for any $G;$
\item $\Bau$ commutes with fltered colimits.
\end{itemize}

The first two statements follow from \cite[Th.4.1]{akhtiamov2021right}. For a finite group $\Bau(G)=1.$ For infinite groups the fact  $|\Bau(G)|\leq |G|$ follows from the equation $|\Bau(F(S))| = |F(S)| = |S|$ for an infinite set $S$  \cite[Th.36.2]{baumslag1960some} and the fact that the epimorphism $F(G)\epi G$ induces an epimorphism $\Bau(F(G))\epi \Bau(G).$ The commuting with filtered colimits follows  from the obvious fact that a filtered colimit of rational groups is rational and the formula ${\sf Hom}({\rm colim}\ G_i,H)= {\rm lim}\ {\sf Hom}(G_i,H).$

\subsection{Baumslag rationalization of a free group}

Baumslag did not consider the functor but he studied rational groups, he called them $\mathcal D$-groups \cite{baumslag1960some}, \cite{baumslag1965free}, \cite{baumslag1968residual}. In particular, he studied the free rational group $\Bau(F).$  He proved that for a free group $F$ there is an increasing family of subgroups $(G_\alpha)_{\alpha\leq \alpha_0}$ of $\Bau(F)$ indexed by ordinals $\leq \alpha_0$ for some ordinal $\alpha_0$
$$ F=G_0 \leq G_1 \leq \dots \leq G_\alpha \leq \dots \leq G_{\alpha_0}=\Bau (F)$$
together with a family of elements $(g_\alpha\in G_\alpha)_{\alpha<\alpha_0}$ such that 
\begin{itemize}
\item $G_\alpha= \bigcup_{\beta<\alpha} G_\beta$ for any limit ordinal $\alpha\leq\alpha_0;$
\item for $\alpha<\alpha_0$ the centraliser $A_\alpha=C(g_\alpha,G_\alpha)$ is isomorphic to a subgroup of $\QQ$ and $G_{\alpha+1}$ is isomorphic to the free product with amalgamation
$G_{\alpha+1}\cong G_\alpha *_{A_\alpha} \QQ.$
$$
\begin{tikzcd}
A_\alpha \arrow[r,rightarrowtail] \arrow[d,hookrightarrow]
\arrow[dr, phantom, "\ulcorner", very near end]
& \QQ \arrow[d] \\
G_\alpha \arrow[r,hookrightarrow] & G_{\alpha+1}
\end{tikzcd}
$$
\end{itemize}
(see {\cite[Th. 36.1]{baumslag1960some}}).

Baumslag also computed the abelianization of this group. He proved that if $F=F(S)$ is a free group freely generated by a set $S,$ then 
$$\Bau(F)_{\sf ab} \cong \QQ^{\oplus S} \oplus (\QQ/\ZZ)^{\oplus T}$$
where the torsion free summand $\QQ^{\oplus S}$ is generated by the images of $S$ and $T$ is a set, which is countable, if $1<|S|<\infty,$ and $|T|=|S|,$ if $S$ is infinite.
(see \cite[Th. 37.3]{baumslag1960some}). In particular, this implies that for a free group $F$ there is an isomorphism
\begin{equation}\label{eq:rational_free_iso_hom}
H_*(\Bau(F),\QQ) \cong H_*(F,\QQ).    
\end{equation}

\subsection{Equivariantly rational modules} Let $G$ be a group and $M$ be a $G$-module (i.e. $\ZZ[G]$-module).  We say that a $G$-module $M$ is {\it equivariantly rational} if the multiplication map 
\[\cdot (1+g+g^2+\dots+g^{n-1}):M\longrightarrow M\]
is a bijection for any $g\in G$ and $n\geq 1.$ It is easy to see that a equivariantly rational $G$-module is a $\QQ$-vector space. If $G$ is rational, a $G$-module $M$ is equivariantly rational if and only if the semidirect product $G\ltimes M$ is a rational group. 

There is an obvious functor of rationalization on the category of modules ${\sf Mod}(G)$ which is just the ring theoretical localization with respect to the set of elements in $\ZZ[G]:$
\begin{equation}
\Sigma=\{ 1+g+g^2+\dots+g^{n-1}\mid g\in G, n\geq 1\}.
\end{equation}

Casacuberta and Peschke generalizaed the equation \eqref{eq:rational_free_iso_hom}. They proved  \cite[Th.8.7, Cor.7.3]{casacuberta1993localizing} that for a free group $F$ and an equivariantly rational $\Bau(F)$-module $M$ the map $F\to \Bau(F)$ induces isomorphisms
\[ H_*(F,M)\cong H_*(\Bau(F),M), \hspace{1cm} H^*(F,M)\cong H^*(\Bau(F),M). \]

\section{Comparison of rationalizations of groups}

Any rational nilpotent group is $H\QQ$-local because it can be obtained from the trivial group as a multiple rational central extension. Hence, the inverse limit of rational nilpotent groups is also $H\QQ$-local. In particular, the $\QQ$-completion of any group $G$ is $H\QQ$-local. From the universal property of $H\QQ$-localization we obtain a homomorphism 
$
\ell_{H\QQ}(G) \longrightarrow \widehat{G}_\QQ
$ (see \eqref{eq:HQ-localization}).  On the other hand, it is easy to see that any $H\QQ$-local group is rational, so we have a map $\Bau(G)\to \ell_{H\QQ}(G).$ Therefore, we obtain two natural transformations 
\begin{equation}
\Bau(G) \longrightarrow \ell_{H\QQ}(G) \longrightarrow \widehat{G}_\QQ.    
\end{equation}

If $H_1(G,\QQ)$ is finite dimensional, then the right hand map is an epimorphism and its kernel is the intersection of the rational lower central series:
\begin{equation}
\frac{\ell_{H\QQ}(G)}{\bigcap_n (\ell_{H\QQ}(G))^\QQ_n  }  \cong \widehat{G}_\QQ.  
\end{equation}

By the definition we have that the map $H_2(G,\QQ)\to H_2(\ell_{H\QQ}(G),\QQ)$ is an epimorphism. In particular, for a free group $F$ we have $H_2(\ell_{H\QQ}(F),\QQ)=0.$ On the other hand $H_2(\widehat{F}_\QQ,\QQ)\ne 0$ for a free group of rank $\geq 2.$ Therefore $\ell_{H\QQ}(F)$ and $\widehat{F}_\QQ$ are not isomorphic
$
\ell_{H\QQ}(F) \not\cong \widehat{F}_\QQ.
$ 

On the other hand we know that for a countable group $G$ the group $\Bau(G)$ is at most countable. But for a finitely generated free group $F$ of rank $\leq 2$ the groups   $\widehat{F}_\QQ$ and $\ell_{H\QQ}(F)$ are uncountable.  Therefore, for a finitely generated free group $F$ the groups 
\begin{equation}
    \Bau(F),\ \ell_{H\QQ}(F),\ \widehat{F}_\QQ, \ \ \text{are not isomorphic}.
\end{equation}

\section{Rationalization of nilpotent spaces}

Let $G$ be a group and $M$ be a $G$-module. The module $M$ is called nilpotent if $MI^n=0$ for some $n,$ where $I$ is the augmentation ideal of $\ZZ[G].$ Note that if  $G$ is a nilpotent group, then $M$ is a nilpotent $G$-module if and only if $G\ltimes M$ is a nilpotent group.  

A connected space $X$ is called {\it nilpotent} if $\pi_1(X)$ is a nilpotent group and $\pi_n(X)$ is a nilpotent $\pi_1(X)$-module for any $n$. In other words a space $X$ is nilpotent if and only if $\pi_1(X)\ltimes \pi_n(X)$ is a nilpotent group for any $n.$ 

Denote by $\WW_{H\QQ}$ the class of all rational homology equivalences i.e. the class of all maps $f:X\to Y$ that induce an isomorphism $H_*(X,\QQ)\cong H_*(Y,\QQ).$ A nilpotent space $X$ is called {\it rational} if one of the following equivalent conditions hold:
\begin{itemize}
    \item $\tilde H_*(X,\ZZ)$ is rational;
    \item $\pi_*(X)$ is rational;
    \item $X$ is $\WW_{H\QQ}$-local 
\end{itemize}
(see \cite[Ch.V,Prop.3.3]{bousfield1972homotopy}). A map between nilpotent spaces $f: X\to Y$ is called {\it rational homotopy equivalence} if one of the following equivalent conditions hold:
\begin{itemize}
\item $f$ is a rational homology equivalence;
\item $\pi_*(X)\otimes \QQ \to  \pi_*(Y)\otimes \QQ$ is iso.
\end{itemize}

Then a rationalization of a space $X$ is a rational homotopy equivalence to a rational space
\begin{equation}
X\longrightarrow X_\QQ.
\end{equation}
In particular 
\begin{equation}
H_*(X_\QQ,\ZZ)\cong H_*(X,\QQ),  \hspace{1cm} \pi_*(X_\QQ)\cong \pi_*(X)\otimes \QQ.     
\end{equation}

Rationalization exists for any nilpotent space \cite[Th.3A]{hilton2016localization},\cite[Ch.V,Prop.4.2]{bousfield1972homotopy}, \cite{bousfield1976pl}. It is easy to see that it is the $\WW_{H\QQ}^{\sf nil}$-localization in the homotopy category of nilpotent spaces, where $\WW_{H\QQ}^{\sf nil}$ is the class of rational homotopy equivalences. Therefore the map $X\to X_\QQ$  satisfies two universal properties: (1) it is the universal map (initial) to a rational space; (2) it is the universal (terminal) rational homotopy equivalence. 

\section{Bousfield-Kan $\QQ$-completion of spaces}

Here by a space we mean a simplicial set. 

For a set $S$ we denote by $\QQ^{(S)}$ the $\QQ$-vector space freely generated by elements of $S.$ For a pointed set $S=(S,s_0)$  we denote by $\QQ(S)$ the quotient $\QQ^{(S)}/(\QQ \cdot s_0).$ This defines a functor on the category of pointed sets
\begin{equation}
    \QQ : {\sf Sets}_* \longrightarrow {\sf Sets}_*.
\end{equation}
This functor has a natural structure of a monad with the obvious unit $\eta_S: S \to \QQ(S)$ and multiplication $\mu:\QQ(\QQ(S))\to \QQ(S),$ whose restriction on the basis $\QQ(S)\!\setminus\! \{0\}$ of $\QQ(\QQ(S))$ is identical. Then any pointed space $S$ has a canonical cosimplicial resolution corresponding to this monad that we denote by $\underline \QQ (S).$

For a pointed simplicial set $X$ we can apply this construction component-wise and obtain a cosimplicial space $\underline \QQ (X).$ Then the Bousfield-Kan $\QQ$-completion $\QQ_\infty X$ is the total space of the cosimplicial space $\underline \QQ (X).$
\begin{equation}
\QQ_\infty X = {\sf Tot}(\underline \QQ (X))    
\end{equation}
(see \cite{bousfield1972homotopy} for details). 

Equivalently Bousfield-Kan $\QQ$-completion of connected spaces can be defined via simplicial groups. The homotopy category of connected pointed spaces is equivalent to the homotopy category of simplicial groups 
\begin{equation}
(\text{connected pointed spaces}) \simeq (\text{simplicial groups}).    
\end{equation}
The equivalence is given by the functors of Kan loop group $X\mapsto \mathcal{G}(X)$ and the functor of simplicial classifying space of a simplicial group $G\mapsto \overline{\mathcal{W}}(G).$
Then the Bousfield-Kan $\QQ$-completion a be defined as the composition of these equivalences with the functor of component-wise $\QQ$-completion on simplicial groups \cite[Ch.IV,\S 4]{bousfield1972homotopy}.
\begin{equation}
\begin{tikzcd}
(\text{connected pointed spaces})
\arrow[rr,"\QQ_\infty"] \ar[d,"\mathcal G"]
&& (\text{connected pointed spaces}) \\ 
\text{(simplicial groups)} \arrow[rr,"\widehat{(-)}_\QQ"]
&& \text{(simplicial groups)} \ar[u,"\overline{\mathcal{W}}"]
\end{tikzcd}
\end{equation}

The $\QQ$-completion can be also defined ``axiomatically'' via towers of fibrations \cite[\S  III.6.3 ]{bousfield1972homotopy}, \cite[\S  12.1]{bousfield1976pl}. Let 
\begin{equation}
\begin{tikzcd}
X \ar[d]  & X \ar[l,"{\sf id}"'] \ar[d] & X \ar[l,"{\sf id}"'] \ar[d] & X \ar[l,"{\sf id}"'] \ar[d] & \dots \ar[l,"{\sf id}"']  \\
N_1 & N_2 \ar[l,twoheadrightarrow] & N_3 \ar[l,twoheadrightarrow] & N_4 \ar[l,twoheadrightarrow] & \dots \ar[l,twoheadrightarrow]
\end{tikzcd}    
\end{equation}
be a commutative diagram such that
\begin{enumerate}
\item $N_i$ are $\QQ$-nilpotent spaces; 
\item the map $\varinjlim H^*(N_i,\QQ)\to H^*(X,\QQ)$ is an isomorphsm; 
\item $N_i\to N_{i-1}$ is a fibration.
\end{enumerate}
Then the map $X\to \QQ_\infty X$ is isomorphic to the map $X\to \varprojlim N_i.$

The Bousfield-Kan $\QQ$-completion can be used to detect rational homology equivalences. For a map $f:X\to Y$ the following statements are equivalent: 
\begin{itemize}
    \item $f$ is a rational homology equivalence;
    \item $\QQ_\infty f : \QQ_\infty X \to \QQ_\infty Y$ is a weak homotopy equivalence. 
\end{itemize}

The $\QQ$-completion of a nilpotent space is the ordinary  rationalization but it is not idempotent for general spaces. A space is called $\QQ$-good if the map $X\to \QQ_\infty X$ is a rational homology equivalence and $\QQ$-bad otherwise. For $\QQ$-bad spaces $\QQ_\infty (\QQ_\infty X)$ is not homotopy equivalent to $\QQ_\infty(X)$ \cite[Ch.I,5.1]{bousfield1972homotopy} It is known that the wedge of two circles $S^1\vee S^1$ is $\QQ$-bad \cite{ivanov2019finite}. Hence 
\begin{equation}
 \QQ_\infty(\QQ_\infty(S^1\vee S^1)) \not \simeq \QQ_\infty(S^1\vee S^1).
\end{equation}

For any space $X$ there is an epimorphism 
\begin{equation}
    \pi_1(\QQ_\infty X) \epi \widehat{ \pi_1(X)}_\QQ
\end{equation}
which is generally not an isomorphism (compare with \cite{ivanov2021right}).

\section{Sullivan's rationalization of spaces}

The strength of the classical rational homotopy theory is due to linear algebraic models of spaces, in particular, minimal Sullivan models \cite{felix2012rational},  \cite{halperin2015rational}. A rationalization of a space can be defined on this language. Namely, Sullivan's rationalization of a space $X$ is the geometric realization of the minimal Sullivan model of $X.$  For spaces of finite rational type Sullivan's rationalization coincides with Bousfield-Kan's $\QQ$-completion.  

Let us explain it in more details following Bousfield and Gugenheim \cite{bousfield1976pl} (see also \cite{halperin2015rational}). By a graded algebra  we mean a non-nenegatively graded algebra $A=\bigoplus_{n\geq 0} A^n.$ A differential graded algebra is a graded algebra with a differential $\partial:A\to A$ of degree $1$ (i.e. a linear map such that $\partial(A^n)\subseteq A^{n+1};$ $\partial^2=0$ and $\partial(ab)=\partial(a)b+(-1)^{|a|}a \partial(b)$ for all homogeneous $a,b\in A$). 
A graded algebra $A$ is called commutative if it satisfies $ab=(-1)^{|a||b|}ba$ for any homogeneous $a,b\in A.$  Cdg-algebra is an abbreviation for commutative differential graded algebra. A morphism of cdg-algebras $f:A\to B$ is called quasi-isomorphism, if it induces an isomorphism on cohomology $H^*(A)\cong H^*(B).$

\subsection{Minimal cdg-algebras}
For any cdg-algebra $A$ and $n\geq 0$ we denote by $A(n)\subseteq A$ the subalgebra generated by $A^0,A^1,\dots,A^n, \partial(A^n);$ and let $A(-1)=\QQ\cdot 1_A \subseteq A_0.$ So we obtain a filtration by cdg-subalgebras
\begin{equation}
\QQ\cdot 1_A = A(-1)\subseteq A(0) \subseteq A(1)\subseteq  \dots     \subseteq A.
\end{equation}
Further, for any $m\geq 0$ we consider a cdg-subalgebra $A(n,m)\subseteq A(n)$ defined inductively: $A(n,0)=A(n-1)$ and $A(n,m+1)$ is the subalgebra generated by $A(n,m)$ and the set $\{ a\in A^n\mid \partial(a)\in A(n,m)\}$
\begin{equation}
A(n-1)=A(n,0)\subseteq A(n,1) \subseteq A(n,2) \subseteq  \dots \subseteq A(n).   
\end{equation}

A cdg-algebra $M$ is called {\it minimal} if
\begin{itemize}
\item it is connected $\QQ\cong M^0;$
\item it is a free commutative graded algebra\\ (generated by some positively graded vector space);
\item $M(n)=\bigcup_{m} M(n,m)$ for any $n\geq 1.$
\end{itemize}
If a cdg-algebra $M$ is simply connected (i.e. $M^0=\QQ$ and $M^1=0$), then $M$ is minimal if and only if it is free commutative graded algebra and $\partial(M)\subseteq \overline{M}\cdot \overline{ M},$ where $\overline{M}=\bigoplus_{n\geq 1}M^n.$

For any homologically connected (i.e. $H^0(A)=\QQ$) cdg-algebra $A$ there exists a quasi-isomorphism from a minimal cdg-algebra $M\to A.$ Moreover, the minimal cdg-algebra $M$ is defined uniquely up to isomorphism. However, $M$ is not natural in $A$ but  it is natural ``up to homotopy''.  

\begin{remark}
There is a model structure on the category of cdg-algebras, where weak equivalences are quasi-isomorphisms and fibrations are surjective morphisms. Minimal algebras are cofibrant with respect to this model structure. Moreover, any cofibrant cdg-algebra $C$ can be decomposed as a tensor product   
\begin{equation}
C\cong M\otimes \bigotimes_{\alpha} T(n_\alpha),    
\end{equation}
where $M$ is a minimal cdg-algebra and for $n\geq 0$ by  $T(n)$ we denote the cdg-algebra freely generated by two elements $a,\partial a,$ where $|a|=n, |\partial a|=n+1$  and $\partial$ is the unique differential such that $\partial(a)=\partial a$ (see \cite[Proposition 7.11]{bousfield1976pl}).
\end{remark}

\subsection{PL de Rham complex}

For any $n\geq 0$ we denote by $\triangledown_n^*$ the cdg-algebra generated by elements $t_0,\dots,t_n$ of  degree $1$ and elements $\partial t_0,\dots, \partial t_n$ of degree $2$ modulo relations 
\begin{equation}
    t_0+t_1+\dots +t_n=1, \hspace{1cm} \partial t_0 + \partial t_1 +\dots + \partial t_n=0
\end{equation}
with the unique differential $\partial : \triangledown_n^m \to \triangledown_n^{m+1}$ such that $\partial(t_i)=\partial t_i.$ We also consider a simlicial object $\triangledown^*_*$ in the category of cdg-algebras, whose components are $\triangledown_n^*$ and the face and degeneracy maps are defined by 
\begin{equation}
d_j(t_i) = 
\begin{cases} 
t_i, & i<j, \\ 
0, &  i=j,\\
t_{i-1}, & i>j,
\end{cases} 
\hspace{1cm}
s_j(t_i) = 
\begin{cases}
t_i, & i<j, \\
t_i+t_{i+1}, & i=j, \\
t_{i+1}, & i>j.
\end{cases}
\end{equation} 
Note that for each fixed $n$ we have a cdg-algebra $\triangledown_n^*$ and for each fixed $m$ we have a simplicial vector space $\triangledown^m_*.$ Then we can define two adjoint functors 
\begin{equation}\label{eq:A_PL}
A_{PL} : \mathsf{sSets} \leftrightarrows \mathsf{cdga}^{op}: F
\end{equation}
such that 
\begin{equation}
A_{PL}(X)^m={\sf Hom}_{\sf sSets}(X,\triangledown_*^m),  \hspace{1cm} F(A)_n = {\sf Hom}_{\sf cdga}(A,\triangledown_n^*),
\end{equation}
where $A_{PL}(X)=\bigoplus_{m\geq 0} A_{PL}(X)^m$ is considered as a cdg-algebra with the pointwise product and $F(A)$ is in fact the composition of two functors $\triangledown:\Delta^{\sf op} \to {\sf cdga}$ and $ {\sf Hom}(A,-):  {\sf cdga}  \to {\sf Sets}.$
The cdg-algebra $A_{PL}(X)$ is called the PL de Rham complex of $X$ or Sullivan de Rham complex of $X.$ The cohomology of $A_{PL}(X)$ is isomorphic to $H^*(X,\QQ).$

\subsection{Sullivan's   rationalization} 

For a connected simplicial set $X$ we consider a quasi-isomorphism $M(X)\to A_{PL}(X)$ from a minimal cdg-algebra $M(X).$ The algebra $M(X)$ is called the Sullivan minimal model of $X.$ $M(X)$ is not natural by $X$ in the category of cdg-algebras but it is natural ``up to homotopy''. Sullivan's rationalization of $X$ is the geometric rationalization of Sullivan's minimal model
\begin{equation}
{\sf Sul}(X) =F(M(X)).    
\end{equation}
The composition of the unit of the adjunction   \eqref{eq:A_PL} 
$\eta : X\to F(A_{PL}(X))$ and the map $F(A_{PL}(X)) \to F(M(X))$ defines the map \begin{equation}
X \longrightarrow {\sf Sul}(X).    
\end{equation}

If $X$ is a space of rational finite type (i.e. $H^n(X,\QQ)$ is finite dimensional for any $n$), then Sullivan's rationalization coincides with Bousfield-Kan' $\QQ$-completion (\cite[Th.12.2]{bousfield1976pl})
\begin{equation}
{\sf Sul}(X)\cong \QQ_\infty X.    
\end{equation}

A similar construction of rationalization can be done via Lie algebra models and it also coincides with  Bousfield-Kan's  $\QQ$-completion \cite{buijs2020lie}.

\section{Homological rationalization of spaces}

Bousfield developed a theory of $h_*$-localization of spaces with respect to any homology theory $h_*$ \cite{bousfield1975localization}. The homology rationalization of a space is the $H_*(-,\QQ)$-localization.

Recall that we denote by $\WW_{H\QQ}$ the class of rational homology equivalences in the homotopy category of spaces i.e. the class of maps $f:X\to Y$ that induce an isomorphism $H_*(X,\QQ)\cong H_*(Y,\QQ).$ Then the homology rationalization of a space is the $\WW_{H\QQ}$-localization. It is denoted by 
\begin{equation}
X \longrightarrow L_{H\QQ}(X).    
\end{equation}
It exists for any space \cite{bousfield1975localization} and defines a functor on the homotopy category.

We say that a space is homologically rational if it is $\WW_{H\QQ}$-local. Hence the homological rationalization satisfies two universal properties: (1) it is the universal rational homology equivalence; (2) it is the universal map to a homologically rational space. 

The homology rationalization corresponds to the $H\QQ$-localization of groups in the following sense. 
The fundamental group of a homological rationalization of a space $X$ is the $H\QQ$-localization of the fundamental group of $X.$
\begin{equation}
\pi_1(L_{H\QQ}(X)) = \ell_{H\QQ}(\pi_1(X)).
\end{equation}

There is a description of connected homologically rational spaces on the language of their homotopy groups. Namely, a connected space $X$ is homologically rational if and only if 
\begin{itemize}
    \item $\pi_n(X)$ is $H\QQ$-local group for $n\geq 1$;
    \item  $\pi_n(X)$ is $H\ZZ$-local as a $\pi_1(X)$-module for $n\geq 2$
\end{itemize}
(see \cite[Th.5.5]{bousfield1975localization}). 

There is a transfinite limit construction for homology rationalization of a space \cite{dror1977long}. 

Note that for the moment there is no a description of $L_{H\QQ}$ in terms of simplicial groups, and the following question is open.

\begin{question}
It the space $L_{H\QQ}(S^1\vee S^1)$ aspherical?
\end{question}

\section{$\Omega$-rationalization of spaces}

\subsection{Definition} Casacuberta and Peschke developed their own theory of localization of spaces \cite{casacuberta1993localizing} (see also \cite{peschke1989localizing}) further developed by  Bastardas and  Casacuberta  \cite{bastardas2001homotopy}, which was initially hinted by Farjoun and Bousfield \cite[Example 7.3]{bousfield1977constructions}.

A connected space $X$ is called {\it $\Omega$-rational} if one of the following equivalent conditions hold:
\begin{itemize}
\item the power map $(-)^m:\Omega X \to \Omega X$ is a homotopy equivalence for any $m\geq 1;$
\item $\pi_1(X)$ is a rational group and $\pi_n(X)$ is an equivariantly rational $\pi_1(X)$-module for $n\geq 2;$
\item $\pi_1(X)$ is  rational and $\pi_1(X)\ltimes \pi_n(X)$ is rational for any $n\geq 2.$
\end{itemize}

Then the $\Omega$-rationalization of a connected space $X$ is the universal map to an $\Omega$-rational space (in the homotopy category)
\begin{equation}
X\longrightarrow L_{\Omega \QQ}(X).  
\end{equation}
It exists for any connected space. 

\subsection{$\Omega$-rationalization as an enriched localization}
The $\Omega$-rationalization is a particular case of the general theory of $f$-localization of spaces \cite{farjoun2006cellular}, \cite{hirschhorn2009model}. In this theory one should consider the category of spaces as a category enriched over itself, it is not enough to think about the homotopy category. If $f:A\to B$ is a map of pointed spaces, a space $X$ is called $f$-local, if the map
\begin{equation}
{\sf Map}_*(B,X) \longrightarrow {\sf Map}_*(A,X)
\end{equation}
is a weak equivalence. The $f$-localization of spaces is the universal map to an $f$-local space. A space is $\Omega$-rational if and only if it is $f$-local, where $f$ is the wedge of the maps $(-)^m:S^1\to S^1$ for $m\geq 1;$ and $\Omega$-rationalization is just the $f$-localization. So the $\Omega$-rationalization shares all properties of $f$-localizations that can be found in \cite{farjoun2006cellular}.  

\subsection{$\Omega$-rationalization as a localization in the homotopy category}
We described $\Omega$-rational spaces as local spaces with respect to the maps $(-)^m:S^1\to S^1$ in the enriched sense. It is also possible to describe them as local spaces in the strict sense with respect to a wider set of maps in the homotopy category, without any enrichment. Consider the sphere with one added point  $S^n_+=S^n\cup {\sf pt}.$ This additional point is considered as the base point of this space. For $n\geq 2$ we also consider the suspension 
\begin{equation}
    S^n_\tau = \Sigma S^{n-1}_+.
\end{equation}
For convenience we set $S^1_\tau:=S^1.$ The space $S^n_\tau$ has the obvious co-$H$-structure because it is a suspension. There is a homotopy equivalence $S^n_\tau \simeq S^n\vee S^1$ but this is not an equivalence of co-$H$-spaces. The space $S^n_\tau$ is interesting because there is a group isomorphism 
\begin{equation}
\pi_1(X)\ltimes \pi_n(X)\cong[S^n_\tau,X],  
\end{equation}
where the structure of a group on $[S^n_\tau,X]$ comes from the co-$H$-structure on $S^n_\tau.$
Then a space $X$ is $\Omega$-rational if and only if it is local with respect to the maps $(-)^n:S_\tau^n\to S_\tau^n, n\geq 1$ in the homotopy category (here the map $(-)^n:S_\tau^n\to S_\tau^n$ is induced by the map $(-)^n:S^1\to S^1$). 

\subsection{$\Omega$-rational equivalences}
We say that a map of connected spaces $f:X\to Y$ is an {\it $\Omega$-equivalence} if $L_{\Omega \QQ}(f)$ is an isomorphism in the homotopy category. There is a description of the class $\Omega$-equivalences. 
For a group $G$ we set 
\begin{equation}
{\bf R}[G] = \ZZ[\Bau(G)][\Sigma^{-1}],    
\end{equation}
where 
\begin{equation}
\Sigma=\{1+g+g^2+\dots+g^{n-1}\mid g\in \Bau(G), n\geq 1\}.    
\end{equation}
Then a map $f:X\to Y$ is an $\Omega$-equivalence if and only if it induces isomorphisms 
\begin{enumerate}
    \item $\Bau(\pi_1(X))\cong \Bau(\pi_1(Y));$
    \item $H^*(Y,\mathcal{A})\cong H^*(X, \mathcal{A} ),$ for any  local system $\mathcal A$ over ${\bf R}[\pi_1(Y)].$ 
\end{enumerate}
The condition (2) is equivalent to the condition 
\begin{itemize}
    \item[(2')]  $H_*(X,{\bf R}[\pi_1(X)])\cong H_*(Y, {\bf R}[\pi_1(Y)]).$ 
\end{itemize}
(see \cite[Th.3.2]{casacuberta1993localizing}). 

\subsection{$\Omega$-rationalization and simplicial groups}
Casacuberta and Bastardas proved that the $\Omega$-rationalization can be defined via simplicial groups as the component-wise Baumslag rationalization \cite{bastardas2001homotopy}. 

\begin{equation}
\begin{tikzcd}
(\text{connected pointed spaces})
\arrow[rr,"L_{\Omega\QQ}"] \ar[d,"\mathcal G"]
&& (\text{connected pointed spaces}) \\ 
\text{(simplicial groups)} \arrow[rr,"\Bau"]
&& \text{(simplicial groups)} \ar[u,"\overline{\mathcal{W}}"]
\end{tikzcd}
\end{equation}

\section{$\pi_1$-fiberwise rationalization of spaces}

G\'{o}mez-Tato, Halperin and Tanr\'{e} developed a version of rationalizations of spaces that I call it {$\pi_1$-fiberwise rationalization} \cite{gomez2000rational}. They also developed a theory of algebraic models that extends the theory of minimal Sullivan modules. 

A space $X$ is called $\pi_1$-fiberwise rational if its universal cover $\tilde X$ is rational in the classical sense. In other words, a space $X$ is $\pi_1$-fiberwise rational if $\pi_n(X)$ is a $\QQ$-vector space for $n
\geq 2$ ($\pi_1(X)$ is arbitrary, the action of $\pi_1(X)$ on $\pi_n(X)$ is arbitrary). A map $f:X\to Y$ is a $\pi_1$-fiberwise rational homotopy equivalence if it induces an isomorphism $\pi_1(X)\cong \pi_1(Y)$ and isomorphisms $\pi_n(X)\otimes \QQ\cong \pi_n(Y)\otimes \QQ$ for $n\geq 2.$ 

Then the $\pi_1$-fiberwise rationalization is the $\pi_1$-fiberwise rational homotopy equivalence to a $\pi_1$-fiberwise rational space 
\begin{equation}
X\longrightarrow L^{\pi_1}_\QQ(X).
\end{equation}
It exists for any space and can be computed as the fiberwise  $\QQ$-completion applied to the fibration $X\to B\pi_1(X).$

The class of $\pi_1$-fiberwise rational homotopy equivalences can be described in terms of rational chain coalgebra $C_*(X,\QQ)$ and in terms of homology. Namely for a map of connected spaces $f:X\to Y$ the following statements are equivalent \cite[Th.16]{rivera2021rational} (see also \cite{rivera2019functor}):
\begin{itemize}
    \item $f$ is a $\pi_1$-fiberwise rational homotopy equivalence; 
    \item $\pi_1(X)\to \pi_1(Y)$ is an isomorphism and $H_*(X,f^*\mathcal{A})\to H_*(Y,\mathcal{A})$ is an isomorphism for any local system $\mathcal A$ over $\QQ[\pi_1(Y)];$
    \item ${\sf Cobar}(C_*(X,\QQ)) \to {\sf Cobar}(C_*(Y,\QQ))$ is a quasi-isomorphism of dg-algebras.
\end{itemize}

\section{Comparison of rationalizations of spaces}

It is not difficult to check that the Bousfield-Kan $\QQ$-completion is a homologically rational space (see  \cite[Prop. 12.10, Lemma 3.5]{bousfield1975localization} and \cite{dror1977long}). On the other hand, looking on homology description of the classes of equivalences, it is easy to see that there are inclusions:
\begin{equation}
\{\pi_1\text{-fiberwise rational eq.}\}\subseteq \{\Omega\text{-rational eq.} \}\subseteq \{\text{homology rational eq.}\}, 
\end{equation}
and for the classes of spaces:
\begin{equation}
\{\pi_1\text{-fiberwise rational sp.}\}\supseteq \{\Omega\text{-rational sp.} \}\supseteq \{\text{homology rational sp.}\}.    
\end{equation}
This implies that there are natural transformations 
\begin{equation}
L^{\pi_1}_\QQ \longrightarrow L_{\Omega\QQ} \longrightarrow L_{H\QQ} \longrightarrow \QQ_\infty.
\end{equation}

For the wedge of two circles we have
\begin{align}
L^{\pi_1}_\QQ(S^1\vee S^1)&= S^1\vee S^1 \\ 
L_{\Omega\QQ}(S^1\vee S^1)&=K(\Bau(F),1)\\
L_{H\QQ}(S^1\vee S^1)&=?\\
\QQ_\infty(S^1\vee S^1)&=K(\widehat{F}_\QQ,1),
\end{align}
where $F$ is the free group of rank two. 
Indeed, for $L^{\pi_1}_\QQ$ it is obvious. For $L_{\Omega\QQ}, \QQ_\infty $ this follows from the simplicial group construction and the fact that $S^1\vee S^1=\overline{\mathcal W}(F),$ where $F$ is the constant free simplicial group. For $L_{H\QQ}$ it is still an open question, if $L_{H\QQ}(S^1\vee S^1)$ is aspherical or not. However, we know that $\pi_1(L_{H\QQ}(S^1\vee S^1))=\ell_{H\QQ}(F).$ Thus these four rationalizations of $S^1\vee S^1$ are non-homotopy equivalent. 
\section{An example: the classifying space of the Burnside group}
Let $G=B(2,n)$ be the free Burnside group of rank two and exponent $n.$ Since the group is generated by torsion elements, $\Bau(G)=\ell_{HQ}(G)=1.$ Therefore, the spaces $L_{\Omega \QQ}(BG)$ and $L_{H\QQ}(BG)$ are simply connected. However they are not contractible for large enough $n$. Indeed, by the Theorem of Ol'shanskii \cite[Cor.31.2]{ol2012geometry} for large enough $n$ the second homology $H_2(G,\ZZ)$ is a free abelian group of countable rank, and hence,
\begin{equation}
H_2(L_{H\QQ}(BG),\QQ)=H_2(L_{H\QQ(BG)},\QQ)=H_2(G,\QQ) \ne 0.   
\end{equation}
By the Hurewicz theorem we obtain
\begin{equation}
\pi_2(L_{\Omega\QQ}(BG))\ne 0, \hspace{1cm} \pi_2(L_{H\QQ}(BG))\ne 0.
\end{equation}

\printbibliography

\end{document}